\documentclass[12 pt]{amsart}

\usepackage{amscd,amsmath,latexsym,amssymb,amscd}
\usepackage[mathscr]{euscript}
\usepackage[all]{xy}
\usepackage{layout}
\usepackage{pstricks}
\usepackage{enumerate}
\usepackage{dsfont}

\newtheorem{theorem}{Theorem}[section]
\newtheorem{corollary}[theorem]{Corollary}
\newtheorem{lemma}[theorem]{Lemma}
\newtheorem{proposition}[theorem]{Proposition}
\theoremstyle{definition}
\newtheorem{definition}[theorem]{Definition}

\theoremstyle{definition}
\newtheorem{example}[theorem]{Example}
\newtheorem{remark}[theorem]{Remark}

\def\Proof{\medskip\noindent{\bf Proof: }}

\def\Z{\mathbb{Z}}

\def\C{\mathbb{C}}

\def\Q{\mathbb{Q}}
\def\R{\mathbb{R}}
\def\F{\mathbb{F}}
\def\S{\mathbb{S}}
\newcommand{\BONE}{\mathds 1}

\def\Fi{\mathcal{F}}
\def\H{\mathcal{H}}

\def\P{\mathcal{P}}
\def\K{\mathcal{K}}

\def\t{\mathfrak{t}}

\def\m{\mathfrak{m}}

\DeclareMathOperator{\Hom}{Hom}

\def\Ker{\text{Ker}}
\def\Rep{\text{Rep}}

\begin{document}

\title[On the Structure of Spaces of Commuting Elements in
Compact Lie Groups]
{On the Structure of Spaces of Commuting
Elements in Compact Lie Groups}

\author[A.~Adem]{Alejandro Adem$^{*}$}
\address{Department of Mathematics,
University of British Columbia, Vancouver BC V6T 1Z2, Canada}
\email{adem@math.ubc.ca}
\thanks{$^{*}$Partially supported by NSERC}

\author[J.~M.~G\'omez]{Jos\'e Manuel G\'omez}
\address{Department of Mathematics,
Johns Hopkins University, Baltimore, MD 21218, USA}
\email{jgomez@math.jhu.edu}

\begin{abstract}
In this note we study topological invariants of the spaces of
homomorphisms $\Hom(\pi,G)$, where $\pi$ is a finitely generated
abelian group and $G$ is a compact Lie group arising as an arbitrary
finite product of the classical groups 
$SU(r)$, $U(q)$ and $Sp(k)$.
\end{abstract}

\maketitle

\section{Introduction}

Let $\P$ denote the class of compact Lie groups arising as arbitrary
finite products of the classical groups  $SU(r)$, $U(q)$
and $Sp(k)$. In this article we
use methods from algebraic topology to study the spaces of
homomorphisms $\Hom(\pi,G)$ where $\pi$ denotes a finitely generated
abelian group and $G\in\P$. Our main interest is the computation of
invariants associated to these spaces such as their cohomology and
stable homotopy type, as well as their equivariant $K$-theory with
respect to the natural conjugation action. The natural quotient
space under this action is the space of representations 
$\Rep(\pi,G)$, which can be identified with the moduli space of 
isomorphism classes of flat connections on principal $G$--bundles 
over $M$, where $M$ is a compact connected manifold with 
$\pi_1(M) = \pi$. Thus our results provide insight into these 
geometric invariants in the important case when $\pi_1(M)$ is a 
finitely generated abelian group.

\medskip

Our starting point is the observation (see \cite{AG}) that when
$G\in\P$  and $\pi$ is a finitely generated abelian group, the
conjugation action of $G$ on the space of homomorphisms
$\Hom(\pi,G)$ satisfies the following property: for every element
$x\in \Hom(\pi,G)$ the isotropy subgroup $G_{x}$ is connected and of
maximal rank. This property plays a central part in our analysis.
Indeed, let $T\subset G$ be a maximal torus; in general if a compact
Lie group $G$ acts on a compact space $X$ with connected maximal
rank isotropy subgroups then there is an associated action of $W$ on
the fixed--point set $X^{T}$ and many properties of the space $X$
are determined by the action of $W$ on $X^{T}$ (see \cite{AG},
\cite{Hauschild}). For our examples this means that a detailed
understanding of the $W$-action on the subspace
$\Hom(\pi,G)^{T}=\Hom(\pi,T)$ can be used to describe key
homotopy--theoretic invariants for the original space of
homomorphisms.

\medskip

This approach can be used for example to obtain an explicit description of
the number of path--connected components in $\Hom(\pi,G)$. Indeed
we show that if
$\pi=\Z^{n}\oplus A$,
where $A$ is a finite abelian group,
then the number of path--connected components
in $\Hom(\pi,G)$ equals the number of distinct orbits
for the action of $W$ on $\Hom(A,T)$

In \cite{AC} a stable splitting for the spaces of commuting
$n$-tuples in $G$, $\Hom(\Z^{n},G)$, was derived for any Lie group
$G$ that is a closed subgroup of $GL_{n}(\C)$. Here we show that
this splitting can be generalized to the spaces of homomorphisms
$\Hom(\pi,G)$ when $G\in \P$ and $\pi$ is any finitely generated
abelian group. This is done by constructing a  stable splitting  on
$\Hom(\pi,G)^{T}=\Hom(\pi,T)$ and proving that this splitting lifts
to the space $\Hom(\pi,G)$. Suppose that
$\pi=\Z/(q_{1})\oplus\cdots\oplus \Z/(q_{n})$, where $n\ge 0$ and
$q_{1},\dots,q_{n}$ are integers. Here we allow some of the
$q_{i}$'s to be $0$ and in that case $\Z/(0)=\Z$. Thus $\Hom(\pi,G)$
can be seen as the subspace of $G^{n}$ consisting of those commuting
$n$-tuples $(x_{1},\dots,x_{n})$ such that $x_{i}^{q_{i}}=1_{G}$ for
all $1\le i\le n$. For $1\le r\le n$ let $J_{n,r}$ denote the set of
all sequences of the form $\m:=\{1\le m_{1}<\cdots<m_{r}\le n\}$.
Given such a sequence $\m$ let
$P_{\m}(\pi):=\Z/(q_{m_{1}})\oplus\cdots\oplus\Z/(q_{m_{r}})$ be a
quotient of $\pi$. Let $S_{1}(P_{\m}(\pi),G)$ be the subspace of
$\Hom(P_{\m}(\pi),G)$ consisting of those $r$-tuples
$(x_{m_{1}},\dots, x_{m_{r}})$ in  $\Hom(P_{\m}(\pi),G)$ for which at
least one of the $x_{m_{i}}$'s is equal to $1_{G}$.

\begin{theorem}\label{splittings intro}
Suppose that $G\in\P$ and that $\pi$ is a finitely generated
abelian group. Then there is a
$G$-equivariant homotopy equivalence
\[
\Theta:\Sigma \Hom(\pi,G)\to \bigvee_{1\le r\le n}\Sigma
\left(\bigvee_{\m\in J_{n,r}}
\Hom(P_{\m}(\pi),G)/S_{1}(P_{\m}(\pi),G)\right).
\]
\end{theorem}

In Section \ref{stable splittings} we determine the homotopy
type of the stable factors appearing in the previous theorem
for certain particular cases. In particular we determine the
stable homotopy type of $\Hom(\pi,SU(2))$ for any
finitely generated abelian group.

\medskip

Suppose now that $G$ is any compact Lie group.
The fundamental group of the spaces of
homomorphisms of the form $\Hom(\Z^{n},G)$ was computed
in \cite{GPS}.
Let $\BONE\in \Hom(\Z^{n},G)$ be the trivial representation.
If $\BONE$ is chosen as the base point, then by
\cite[Theorem1.1]{GPS} there is a natural isomorphism
$\pi_{1}(\Hom(\Z^{n},G))\cong (\pi_{1}(G))^{n}$.
Here we show that the methods applied in \cite{GPS} can be
used to compute $\pi_{1}(\Hom(\pi,G))$ for any choice of base point
if we further require that $G\in\P$ and that $\pi$ is a finitely
generated abelian group.  Write $\pi$ in the form
$\pi=\Z^{n}\oplus A$, with $A$ a finite abelian group.
Then the space of homomorphisms $\Hom(\pi,G)$ can
naturally be identified as a subspace of the product
$\Hom(\Z^{n},G)\times \Hom(A,G)$. Given
$f\in \Hom(A,T)$ let
\[
\BONE_{f}:=\BONE\times f\in \Hom(\pi,G)\subset \Hom(\Z^{n},G)\times
\Hom(A,G).
\]
Every path--connected component in
$\Hom(\pi,G)$ contains some $\BONE_{f}$ and thus it
suffices to consider the elements of the form $\BONE_{f}$
as base points in $\Hom(\pi,G)$.  With this in mind we have the following.

\begin{theorem}
Let $\pi=\Z^{n}\oplus A$, with $A$ a finite abelian group
and let  $G\in\P$. Suppose $f\in \Hom(A,T)$  and take $\BONE_{f}$ 
as the base point of $\Hom(\pi,G)$.  Then there is a natural 
 isomorphism $\pi_{1}(\Hom(\pi,G))\cong (\pi_{1}(G_{f}))^{n}$
where $G_{f}=Z_{G}(f)$ is the subgroup of elements in $G$ commuting
with $f(x)$ for all $x\in A$.
\end{theorem}

In Section \ref{equivariant K-theory} we study the equivariant
$K$-theory of the spaces of homomorphisms $\Hom(\pi,G)$ with respect
to the conjugation action by $G$. When $\pi$ is a finite group, then
$\Hom(\pi,G)$ is the disjoint union of homogeneous spaces of the
form $G/H$ where $H$ is a maximal rank subgroup. Using this it is
easy to see that $K^{*}_{G}(\Hom(\pi,G))$ is a free module over the
representation ring of rank $|\Hom(\pi,T)|$. This result can be
generalized for finitely generated abelian groups of rank $1$ in the
following way.

\begin{theorem}
Suppose that $G\in \P$ is simply connected and of rank $r$.
Let $\pi=\Z\oplus A$ where $A$ is a finite abelian group.
Then $K_{G}^{*}(\Hom(\pi,G))$ is a free $R(G)$-module of rank
$2^{r}\cdot |\Hom(A,T)|$.
\end{theorem}

It turns out that $K^{*}_{G}(\Hom(\pi,G))$ is not always free as a
module over $R(G)$. In fact, as was pointed out in \cite{AG}, the
$R(SU(2))$-module $K^{*}_{SU(2)}(\Hom(\Z^{2},SU(2)))$ is not free.
However, $K^{*}_{SU(2)}(\Hom(\Z^{2},SU(2)))\otimes \Q$ turns out to
be free as a module over $R(SU(2))\otimes \Q$. The next theorem
shows that a similar result holds for all the spaces of
homomorphisms that we consider here.

\begin{theorem}
Suppose that $G\in \P$ is of rank $r$ and that $\pi$ is a finitely
generated abelian group written in the form $\pi=\Z^{n}\oplus A$,
where $A$ is a finite abelian group. Then
$K_{G}^{*}(\Hom(\pi,G))\otimes \Q$ is a free module over
$R(G)\otimes \Q$ of rank $2^{nr}\cdot |\Hom(A,T)|$.
\end{theorem}

The layout of this article is as follows. In Section \ref{spaces of comm}
some general properties of the spaces of homomorphisms $\Hom(\pi,G)$
are determined. In Section \ref{max rank} we study the cohomology groups
with rational coefficients of these spaces. In Section \ref{stable splittings}
Theorem \ref{splittings intro} is proved and some explicit examples
are computed. In Section \ref{fundamental group} the fundamental
group of the spaces $\Hom(\pi,G)$ are computed for any choice of
base point. Finally, in Section \ref{equivariant K-theory} we study the
problem of computing $K_{G}^{*}(\Hom(\pi,G))$, where $G$ acts
by conjugation on $\Hom(\pi,G)$.

Both authors would like to thank the Centro di Ricerca Matematica
Ennio De Giorgi at the Scuola Normale Superiore in Pisa for inviting
them to participate in the program on Configuration Spaces:
Geometry, Combinatorics and Topology during the spring of 2010.

\section{Preliminaries on spaces of commuting elements}\label{spaces of comm}

Let $\pi$ be a finitely generated discrete group and $G$ a Lie
group. Consider the set of homomorphisms from $\pi$ to $G$,
$\Hom(\pi,G)$. This set can be given a  topology as a subspace of a
finite product of copies of $G$ in the following way. Fix a set of
generators $e_{1},\dots, e_{n}$ of $\pi$ and let $F_{n}$ be the free
group on $n$-letters. By mapping the generators of $F_{n}$ onto the
different $e_{i}$'s we obtain a surjective homomorphism $F_{n}\to
\pi$. This surjection induces an inclusion of sets
$\Hom(\pi,G)\hookrightarrow \Hom(F_{n},G)\cong G^{n}$. This way
$\Hom(\pi,G)$ can be given the subspace topology. It is easy to see
that this topology is independent of the generators chosen for
$\pi$. In case $\pi$ happens to be abelian, then any map $F_{n}\to
\pi$ factors through $F_{n}\to \Z^{n}\to  \pi$ yielding an inclusion
of spaces $\Hom(\pi,G)\hookrightarrow\Hom(\Z^{n},G)\hookrightarrow
G^{n}$. Thus the space of homomorphisms $\Hom(\pi,G)$ can be seen as
a subspace of the space of commuting $n$-tuples in $G$,
$\Hom(\Z^{n},G)$.

\medskip

In this note we collect some facts about these spaces of homomorphisms
in the particular case that $\pi$ is a finitely generated abelian group
and $G$ belongs to a suitable family of Lie groups.   We are
mainly interested in the following family of Lie groups.

\begin{definition}
Let $\P$  denote the collection of all compact Lie groups arising as
finite cartesian products of the groups $SU(r)$, $U(q)$ and $Sp(k)$.
\end{definition}

Whenever $G$ belongs to the family $\P$ the space of
homomorphisms $\Hom(\pi,G)$ satisfies the following crucial
condition as we prove below in Proposition
\ref{max rank for spaces of hom}.

\begin{definition}
Let $X$ be a $G$-space. The action of $G$ on $X$ is said to have
connected maximal rank isotropy subgroups if for every $x\in X$,
the isotropy group $G_{x}$ is a connected subgroup of maximal rank;
that is, for every $x\in X$ we can find a maximal torus $T_{x}$ in $G$
such that $T_{x}\subset G_{x}$.
\end{definition}

\begin{proposition}\label{max rank for spaces of hom}
Suppose that $\pi$ is a finitely generated abelian group and $G\in \P$.
Then the conjugation action of $G$ on $\Hom(\pi,G)$ has connected
maximal rank isotropy subgroups.
\end{proposition}
\Proof Choose generators $e_{1},\dots,e_{n}$ of $\pi$. As pointed
out above we can use these generators to obtain an inclusion of
$G$-spaces $\Hom(\pi,G)\hookrightarrow \Hom(\Z^{n},G)$. Given this
inclusion it suffices to show that the conjugation action of $G$ on
$\Hom(\Z^{n},G)$ has connected maximal rank isotropy groups. In
\cite[Example 2.4]{AG} it was proven that the action of 
$G$ on $\Hom(\Z^{n},G)$ has connected maximal rank isotropy subgroups 
if and only if $\Hom(\Z^{n+1},G)$ is
path--connected. The proposition follows by noting that
$\Hom(\Z^{k},G)$ is path--connected for all $k\ge 0$ whenever $G\in
\P$. \qed

\medskip

Suppose that a compact Lie group $G$ acts on a space $X$ with
connected maximal rank isotropy subgroups.  Choose a maximal
torus $T$ in  $G$ and let $W$ be the Weyl group.  By passing
to the level of $T$-fixed points, the action of $G$ on $X$ induces an
action of the Weyl group $W$ on $X^{T}$.  Many properties of the
action of  $G$ on $X$ are determined  by the action of $W$ on $X^{T}$
as explained in  \cite{Hauschild} and in some situations the
former is completely determined by the latter up to isomorphism.
For example, we can use this approach to produce $G$-CW complex
structures on the spaces of homomorphisms as is proved next.

\begin{corollary}\label{CW complex}
Suppose that $\pi$ is a finitely generated abelian group and $G\in \P$.
Then $\Hom(\pi,G)$ with the conjugation action has the structure of a
$G$-CW complex.
\end{corollary}
\Proof
Since $\pi$ is a finitely generated abelian group it can be
written in the form $\pi=\Z^{n}\oplus A$, where $A$ is a
finite abelian group.
Let $X:=\Hom(\pi,G)$ with the conjugation action of $G$. Note that
$X^{T}=\Hom(\pi,G)^{T}=T^{n}\times \Hom(A,T)$.
Since $\Hom(A,T)$ is a discrete set, it follows that $X^{T}$ has the
structure of a smooth manifold on which $W$ acts smoothly. In particular,
by \cite[Theorem 1]{Illman} it follows that $X^{T}$ has the structure
of a $W$-CW complex. Since the conjugation action of $G$ on
$X$ has connected maximal rank isotropy subgroups then by
\cite[Theorem 2.2]{AG} it follows that this $W$-CW complex
structure on $X^{T}$ induces a $G$-CW complex on $X$.
\qed

\medskip

This approach can also be used to determine explicitly the structure
of these spaces of homomorphisms whenever $\pi$ is a finite abelian
group.

\begin{proposition}\label{finite groups}
Suppose that $\pi$ is a finite abelian group and $G\in \P$. Then
there is a $G$-equivariant homeomorphism
\[
\Phi:\Hom(\pi,G)\to \bigsqcup_{[f]\in \Hom(\pi,T)/W}G/G_{f}.
\]
Here $[f]$ runs through a system of representatives of the
$W$-orbits in $\Hom(\pi,T)$ and each $G_{f}$ is a maximal rank
subgroup with $W(G_{f})=W_{f}$.
\end{proposition}
\Proof
Consider the $G$-space $X:=\Hom(\pi,G)$. Note that
$X^{T}=\Hom(\pi,T)$
is a discrete set endowed with an action of $W$. By decomposing
$X^{T}$ into the different $W$-orbits we obtain a $W$-equivariant
homeomorphism
\[
X^{T}\cong \bigsqcup_{[f]\in \Hom(\pi,T)/W}W/W_{f}.
\]
Here $[f]$ runs through a set of representatives for the action of
$W$ on $\Hom(\pi,T)$. For each $f\in \Hom(\pi,T)$ let $G_{f}$ denote
the subgroup of elements in $G$ commuting with $f(x)$ for all $x\in
\pi$. This group is a maximal rank subgroup in $G$ as $T\subset
G_{f}$. Moreover, by \cite[Theorem 1.1]{Hauschild} it follows that
$W(G_{f})=W_{f}$. Also note that if we let $G$ act on the left on
the homogeneous space $G/G_{f}$ then $(G/G_{f})^{T}=W/W_{f}$. Let
\[
Y=\bigsqcup_{[f]\in \Hom(\pi,T)/W}G/G_{f}.
\]
The left action of $G$  on  $Y$ has  maximal rank isotropy and there
is a $W$-equivariant homeomorphism $\phi:X^{T}\to Y^{T}$. By
\cite[Theorem 2.1]{Hauschild} there is a unique $G$-equivariant
extension $\Phi:X\to Y$  of $\phi$ and this map is in fact a
homeomorphism. \qed

\section{Rational cohomology and path--connected components}\label{max rank}

In this section we explore the set of path connected components and
the rational cohomology groups of the spaces of homomorphisms
$\Hom(\pi,G)$.

\medskip

Suppose that $G$ is a compact connected Lie group and let
$T$ be a maximal torus in $G$.  Assume that $G$ acts on
a space $X$ of the homotopy type of a $G$-CW complex
with maximal rank isotropy subgroups. Consider the
continuous map
\begin{align*}
\phi:G\times X^{T}&\to X \\
(g,x)&\mapsto gx.
\end{align*}

Since $G$ acts on $X$ with maximal rank isotropy subgroups
for every $x\in X$ we can find a maximal torus $T_{x}$
 in $G$ such that $T_{x}\subset G_{x}$. As every pair of
maximal tori in $G$ are conjugate it follows that for
every $x\in X$ we can find some $g\in G$ such that
$gx\in X^{T}$. This shows that $\phi$ is a surjective map.
The normalizer of $T$ in $G$, $N_{G}(T)$ acts on the right
on $G\times X^{T}$ by $(g,x)\cdot n= (gn,n^{-1}x)$ and the map
$\phi$ is invariant under this action. Thus $\phi$ descends to
a surjective map
\begin{align*}
\varphi:G\times_{N_{G}(T)}X^{T}=G/T\times_{W}X^{T}&\to X\\
[g,x]&\mapsto gx
\end{align*}
The map $\varphi$ is not injective in general. Indeed, as was proven
in \cite{Baird}, given $x\in X$ there is a homeomorphism
$\varphi^{-1}(x)\cong G_{x}^{0}/N_{G_{x}^{0}}(T)$, where $G_{x}^{0}$
denotes the path--connected component of $G_{x}$ containing the
identity element.  Let $\F$ be a field with characteristic
relatively prime to $|W|$. Then as observed in \cite{Baird} the
space $G_{x}^{0}/N_{G_{x}^{0}}(T)$ has $\F$--acyclic cohomology. The
Vietoris-Begle theorem shows that $\varphi$ induces an isomorphism
in cohomology with $\F$-coefficients. As a consequence we obtain the
following proposition (first proved in \cite{Baird}).

\begin{proposition}\label{rational cohomology}
Suppose that $G$ is a compact connected Lie group
acting on a spaces $X$ with maximal rank isotropy subgroups.
If $\F$ is a field with characteristic relatively prime to $|W|$
then
$H^{*}(X;\F)\cong H^{*}(G/T\times_{W} X^{T};\F)\cong
H^{*}(G/T\times X^{T};\F)^{W}$.
\end{proposition}

\begin{remark}
Suppose that $G$ acts on $X$ with {\em{connected}} maximal rank
isotropy groups. As pointed out above the map $\varphi$ is not
injective in general since $\varphi^{-1}(x)\cong
G_{x}^{0}/N_{G_{x}^{0}}(T)$ for $x\in X$.  Under the given
hypothesis we have $G_{x}^{0}=G_{x}$. By \cite[Theorem
1.1]{Hauschild} the assignment $(H)\mapsto (WH)$ defines  a one to
one correspondence between the set of conjugacy classes of isotropy
subgroups of the action of $G$ on $X$ and the set of conjugacy
classes of isotropy subgroups of the action of $W$ on $X^{T}$. 
Thus the different isotropy
subgroups of the action of $W$ on $X^{T}$ determine how far the map
$\varphi$ is from being injective. In particular, if $W$ acts freely
on $X^{T}$ then $\varphi$ is a continuous bijection and thus a
homeomorphism if for example $X^{T}$ is compact.
\end{remark}

Suppose now that $G\in \P$ and let $\pi$ be a finitely
generated abelian group.  By Proposition
\ref{max rank for spaces of hom} the conjugation action of
$G$ on $\Hom(\pi,G)$ has connected maximal rank isotropy
subgroups. In this case
$\Hom(\pi,G)^{T}=\Hom(\pi,T)$.
As a consequence of the previous result the following
is obtained.
\begin{corollary}\label{cohomology hom}
Suppose that $G\in \P$ and let $\pi$ be a finitely generated
abelian group. Then there is an isomorphism
$H^{*}(\Hom(\pi,G);\Q)\cong H^{*}(G/T\times \Hom(\pi,T);\Q)^{W}$.
\end{corollary}

As an application of Corollary \ref{cohomology hom}
the following can be derived.

\begin{corollary}
Suppose that $G\in \P$ and let $\pi$ be a finitely generated
abelian group written in the form $\pi=\Z^{n}\oplus A$.
Then the number of path--connected components
in $\Hom(\pi,G)$ equals the number of different orbits
of the action of $W$ on $\Hom(A,T)$
\end{corollary}

\section{Stable splittings}\label{stable splittings}

In this section we show that the fat wedge filtration on a
finite product of copies of $G$ induces a natural filtration on the
spaces of homomorphisms $\Hom(\pi,G)$. It turns out that
this filtration splits stably after one suspension whenever
$\pi$ is a finitely generated abelian group and $G\in \P$.

\medskip

Suppose that $\pi$ is a  finitely generated abelian group. Using the
fundamental theorem of finitely generated abelian groups $\pi$ can
be written in the form
\[
\pi=\Z/(q_{1})\oplus\cdots\oplus \Z/(q_{n}),
\]
where $n\ge 0$ and $q_{1},\dots,q_{n}$ are integers. Here we
allow some of the $q_{i}$'s to be $0$ and in that case $\Z/(0)=\Z$.
This way we can see $\Hom(\pi,G)$ as the subspace of $G^{n}$
consisting of those commuting $n$-tuples $(x_{1},\dots,x_{n})$
such that $x_{i}^{q_{i}}=1_{G}$ for all $1\le i\le n$. The
fat wedge filtration on $G^{n}$ induces a
natural filtration on the space of homomorphisms $\Hom(\pi,G)$.
To be more precise, for each $1\le j\le n$ let
\[
S_{j}(\pi,G)=\{(x_{1},\dots,x_{n})\in \Hom(\pi,G)\subset G^{n}
~~|~~ x_{i}=1_{G} \text{ for at least } j \text{ of the } x_{i}\text{'s}\}.
\]
This way we obtain a filtration of $\Hom(\pi,G)$
\begin{equation}\label{filtration1}
\{(1_{G},...,1_{G})\}=S_{n}(\pi,G)\subset
S_{n-1}(\pi,G)\subset\cdots \subset S_{0}(\pi,G)=\Hom(\pi,G).
\end{equation}
Note that each $S_{j}(\pi,G)$ is invariant under the conjugation
action of $G$. In particular each $S_{j}(\pi,G)$ can be seen as a
$G$-space that has connected maximal rank isotropy
subgroups. On the level of the $T$-fixed points the
filtration (\ref{filtration1}) induces a filtration of
$\Hom(\pi,G)^{T}$
\begin{equation}\label{filtration2}
\{(1_{G},...,1_{G})\}=S_{n}(\pi,G)^{T}\subset
S_{n-1}(\pi,G)^{T}\subset\cdots \subset S_{0}(\pi,G)^{T}
=\Hom(\pi,G)^{T}.
\end{equation}
For each $1\le i\le n$ consider
$\Hom(\Z/q_{i},T)=\{t\in T ~~|~~ t^{q_{i}}=1\}$.
Note that each $\Hom(\Z/q_{i},T)$ is a space endowed with the
action of $W$. Whenever $q_{i}=0$ we have $\Hom(\Z/q_{i},T)=T$
and if $q_{i}\ne 0$ then $\Hom(\Z/q_{i},T)$ is a discrete set.
Since $T$ is abelian it follows that
\[
\Hom(\pi,G)^{T}=\Hom(\pi,T)=\Hom(\Z/q_{1},T)\times\cdots
\times \Hom(\Z/q_{n},T).
\]
Moreover, the filtration (\ref{filtration2}) is precisely the
fat wedge filtration of $\Hom(\pi,G)^{T}$  where we identify
$\Hom(\pi,G)^{T}$ with the above product.
It is well known that the fat wedge filtration on a product of
spaces splits stably after one suspension. More precisely, for each
$0\le j\le n-1$ we can find a continuous map
\[
r_{j}:\Sigma S_{j}(\pi,G)^{T}\to \Sigma S_{j+1}(\pi,G)^{T}
\]
in such a way that there is a homotopy $h_{j}$ between
$r_{j} \circ \Sigma(i_{j})$ and $1_{\Sigma(S_{j+1}(\pi,G)^{T})}$.
Here
\[
i_{j}:S_{j+1}(\pi,G)^{T}\to S_{j}(\pi,G)^{T}
\]
denotes the inclusion map. Moreover, both the map $r_{j}$ and the
homotopy $h_{j}$ can be arranged in such a way that they are
$W$-equivariant. The $W$-action that we have in sight is the diagonal
action of $W$ on the product $\Hom(\Z/q_{1},T)\times\cdots
\times \Hom(\Z/q_{n},T)$.  Consider the action of $G$ on
$\Sigma \Hom(\pi,G)$ with $G$ acting trivially on the suspension
component. This action has connected maximal rank isotropy
subgroups and  $(\Sigma \Hom(\pi,G))^{T}=\Sigma \Hom(\pi,T)$.
By \cite[Theorem 2.1]{Hauschild} we can find a unique
$G$-equivariant extension
\[
R_{j}:\Sigma S_{j}(\pi,G)\to \Sigma S_{j+1}(\pi,G)
\]
of $r_{j}$ and a unique $G$-equivariant homotopy $H_{j}$ between
$R_{j} \circ \Sigma(I_{j})$ and
$1_{\Sigma(S_{j+1}(\pi,G))}$ extending $h_{j}$. Here
$I_{j}:S_{j+1}(\pi,G)\to S_{j}(\pi,G)$ as before denotes
the inclusion map.

Let $J_{n,r}$ denote the set of all sequences of the form
$\m:=\{1\le m_{1}<\cdots<m_{r}\le n\}$.
Note that $J_{n,r}$ contains precisely
$\binom{{n}}{{r}}$ elements. Given such a sequence $\m$,
there is an associated abelian group
$P_{\m}(\pi):=\Z/(q_{m_{1}})\oplus\cdots\oplus\Z/(q_{m_{r}})$
obtained as a quotient of $\pi$ and also a
$G$-equivariant projection map
\begin{align*}
P_{\m}:\Hom(\pi,G)&\to \Hom(P_{\m}(\pi),G)\\
(x_{1},...,x_{n})&\mapsto (x_{m_{1}},...,x_{m_{r}}).
\end{align*}

The above can be used to prove the following theorem.

\begin{theorem}\label{stable splitting}
Suppose that $G\in\P$ and that $\pi$ is a finitely generated
abelian group. Then there is a $G$-equivariant homotopy equivalence
\[
\Theta:\Sigma \Hom(\pi,G)\to \bigvee_{1\le r\le n}\Sigma
\left(\bigvee_{\m\in J_{n,r}}
\Hom(P_{\m}(\pi),G)/S_{1}(P_{\m}(\pi),G)\right).
\]
\end{theorem}
\Proof
Note that each $S_{j}(\pi,G)^{T}$ has the homotopy type of a
$W$-CW complex and this implies that each $S_{j}(\pi,G)$
has the homotopy type of a $G$-CW complex by
\cite[Theorem 2.2]{AG}. The different maps $R_{j}$ and the
homotopies $H_{j}$ induce a $G$-equivariant homotopy equivalence
\[
\Sigma\Hom(\pi,G)\simeq \bigvee_{0\le r\le n-1}\Sigma
S_{r}(\pi,G)/S_{r+1}(\pi,G)=\bigvee_{1\le r\le n}\Sigma
S_{n-r}(\pi,G)/S_{n-r+1}(\pi,G).
\]
To finish the theorem we will show that for each $1\le r\le n$
there is a $G$-equivariant homotopy equivalence
\[
S_{n-r}(\pi,G)/S_{n-r+1}(\pi,G)\simeq
\bigvee_{\m\in J_{n,r}}
\Hom(P_{\m}(\pi),G)/S_{1}(P_{\m}(\pi),G).
\]
To see this note that the different projection maps
$\{P_{\m}\}_{\m\in J_{n,r}}$
can be assembled to obtain a $G$-map
\begin{align*}
\eta:\Hom(\pi,G)&\to
\prod_{\m\in J_{n,r}}\Hom(P_{\m}(\pi),G)/S_{1}(P_{\m}(\pi),G)\\
(x_{1},...,x_{n})&\mapsto \{\bar{P}_{\m}
(x_{1},...,x_{n})\}_{\m\in J_{n,r}}.
\end{align*}
The map $\eta$ sends $S_{n-r}(\pi,G)$ onto
$\bigvee_{\m\in J_{n,r}}{\Hom(P_{\m}(\pi),G)/S_{1}(P_{\m}(\pi),G)}$
and $S_{n-r+1}(\pi,G)$ is mapped onto the base point.
It is easy to see that $\eta$ induces a $G$-equivariant homeomorphism
\[
S_{n-r}(\pi,G)/S_{n-r+1}(\pi,G)\cong
\bigvee_{\m\in J_{n,r}}
\Hom(P_{\m}(\pi),G)/S_{1}(P_{\m}(\pi),G)
\]
and the theorem follows.
\qed

\medskip

\noindent{\bf{Remark:}} A case of particular importance in the
previous theorem is $\pi=\Z^{n}$. In this case $\Hom(\Z^{n},G)$ is
precisely the space of commuting ordered $n$-tuples in $G$. The
previous theorem provides a simple proof for the stable equivalence
provided in \cite{AC} for the spaces $\Hom(\Z^{n},G)$ whenever $G\in
\P$.

\begin{example}\label{commuting in SU(2)}
Suppose that $\pi=\Z^{n}$. Let $1\le r\le n$. For any
$\m\in J_{n,r}$ we have
\[
\Hom(P_{\m}(\pi),G)/S_{1}(P_{\m}(\pi),G)\cong \Hom(\Z^{r},G)/S^{1}_{r}(G),
\]
where $S^{1}_{r}(G)\subset \Hom(\Z^{n},G)$ is the subspace of those
commuting $n$-tuples $(x_{1},\dots,x_{n})$ with  at least one
of the $x_{i}$ equal to $1_{G}$. These stable factors were identified
independently in \cite{ACG}, \cite{BJS} and \cite{Crabb}  in the particular
case where $G=SU(2)$. Let
$n\lambda_{2}$ denote the the Whitney sum of $n$-copies of the
canonical vector bundle over $\R P^{2}$ and
let $s_{n}$ denote its zero section. Then
\begin{equation*}
\Hom(\Z^{n},SU(2))/S^{1}_{n}(SU(2))\cong\left\{
\begin{array}{cl}
\S^{3} & \text{if }n=1, \\
(\R P^{2})^{n\lambda_{2}}/s_{n}(\R P^{2})& \text{if }n\ge 2.
\end{array}
\right.
\end{equation*}
\end{example}

\medskip

\begin{example}
Suppose now that
$\pi=\Z/(q_{1})\oplus\cdots\oplus \Z/(q_{n})$
is any finitely generated abelian group and $G=SU(2)$.
Let $T$ be the maximal torus consisting of $2\times 2$
diagonal matrices with entries in $\S^{1}$ and determinant $1$.
In this case $W=\Z/2$ acts by permuting the diagonal entries of
elements in $T$.  Next we determine the stable factors of the
form
\[
\Hom(P_{\m}(\pi),SU(2))/S_{1}(P_{\m}(\pi),SU(2)),
\]
where $\m=\{1\le m_{1}<\cdots<m_{r}\le n\}$ is fixed.
We consider the following cases.

\begin{list}{\labelitemi}{\leftmargin=1em}

\item Suppose $P_{\m}(\pi)$ is a finite group so that $q_{m_{i}}\ne 0$
for all $1\le i\le r$. Assume further that at least
one of the $q_{m_{i}}$'s is odd.
By Proposition \ref{finite groups} there is a homeomorphism
\[
\Hom(P_{\m}(\pi),SU(2))\cong \bigsqcup_{[f]\in
\Hom(P_{\m}(\pi),T)/W}G/G_{f}.
\]
Here $[f]$ runs through all the $W$-orbits in
$\Hom(P_{\m}(\pi),T)$. In this case
\[
G/G_{f}=G/T\cong \S^{2}
\]
for all orbits corresponding to elements $f$ for
which $W_{f}$ is trivial. On the other hand, when $f$ is fixed by
$W$ the corresponding orbit is $G/G_{f}=G/G=*$. Since we
are assuming that one of the $q_{m_{i}}$'s is odd, then
every $f\in\Hom(P_{\m}(\pi),T)$ corresponding to
$r$-tuples $(x_{m_{1}},\dots,x_{m_{r}})$ in $\Hom(P_{\m}(\pi),SU(2))$
with $x_{m_{i}}\ne 1_{G}$ for all $i$ satisfies $W_{f}=1$.
This shows that
\[
\Hom(P_{\m}(\pi),SU(2))/S_{1}(P_{\m}(\pi),SU(2))\cong
(\bigsqcup_{A(\m,\pi)} \S^{2})_{+}.
\]
Here $A(\m,\pi)$ is the number of $W$-orbits in $\Hom(P_{\m}(\pi),T)$
corresponding to $r$-tuples that don't contain
the element $1$. This number is precisely
\[
A(\m,\pi)=\frac{1}{2}(q_{m_{1}}-1)\cdots(q_{m_{r}}-1).
\]

\item Suppose now that $q_{m_{i}}\ne 0$ is even for all $1\le i\le r$.
In this case we have two possibilities for the $W$-orbits in
$\Hom(P_{\m}(\pi),T)$. If $[f]$ represents the orbit $[(-1,\dots,-1)]$
then $W_{f}=W$ and the corresponding orbit is $G/G_{f}=*$.
For all other orbits $[f]\in \Hom(P_{\m}(\pi),T)/W$ corresponding to
$r$-tuples $(x_{m_{1}},\dots,x_{m_{r}})$ in $\Hom(P_{\m}(\pi),SU(2))$
with $x_{m_{i}}\ne 1_{G}$ for all $i$ we have $W_{f}=1$ and
as before $G/G_{f}\cong \S^{2}$. This shows that
\[
\Hom(P_{\m}(\pi),SU(2))/S_{1}(P_{\m}(\pi),SU(2))\cong
(\bigsqcup_{A(\m,\pi)} \S^{2})\sqcup \S^{0},
\]
where now $A(\m,\pi)$ is the number of $W$-orbits
in $\Hom(P_{\m}(\pi),T)$ corresponding
to $r$-tuples in $\Hom(P_{\m}(\pi),T)$ that don't
contain the element $1$ and that are different from
$(-1,\dots,-1)$. This number is precisely
\[
A(\m,\pi)=\frac{1}{2}\left((q_{m_{1}}-1)\cdots(q_{m_{r}}-1)-1\right).
\]

\item We now consider the case where $q_{m_{i}}=0$ for some $1\le i\le r$.
If $q_{m_{i}}=0$ for all $1\le i\le r$ then
\[
\Hom(P_{\m}(\pi),SU(2))/S_{1}(P_{\m}(\pi),SU(2))=
\Hom(\Z^{r},SU(2))/S^{1}_{r}(SU(2))
\]
and these stable factors are as in Example \ref{commuting in SU(2)}.
Suppose then that $q_{m_{i}}\ne 0$ for some $i$. For simplicity
and without loss of generality we may assume that
\[
P_{\m}(\pi)=\Z^{k}\oplus\Z/(q_{m_{k+1}})\oplus\cdots\oplus \Z/(q_{m_{r}})
\]
for some $1\le k <r$ and $q_{m_{i}}\ne 0$ for $k+1\le i\le r$.
Since the inclusion map
$S_{1}(P_{\m}(\pi),SU(2))\hookrightarrow \Hom(P_{\m}(\pi),G) $
is a cofibration, we have
\[
\Hom(P_{\m}(\pi),SU(2))/S_{1}(P_{\m}(\pi),SU(2))\cong
(\Hom(P_{\m}(\pi),SU(2))\setminus S_{1}(P_{\m}(\pi),SU(2)))^{+}.
\]
Here if $X$ is a locally compact space then $X^{+}$ denotes its
one point compactification. Consider the map
\begin{align*}
\varphi_{\m}:G/T\times_{W} \Hom(P_{\m}(\pi),T)&\to \Hom(P_{\m}(\pi),G)\\
(g,(t_{m_{1}},\dots,t_{m_{r}}))&\mapsto
(gt_{m_{1}}g^{-1},\dots,gt_{m_{r}}g^{-1}).
\end{align*}
This map is surjective as the action of $G$ on $\Hom(P_{\m}(\pi),G)$
has connected maximal rank isotropy. Moreover
\[
\varphi_{\m}(g,(t_{m_{1}},\dots,t_{m_{r}}))\in S_{1}(P_{\m}(\pi),SU(2))
\]
if and only if $t_{m_{i}}=1$ for some $1\le i\le r$. Let
$Q(P_{\m}(\pi),T)$ denote the subset of
\[
\Hom(\Z/(q_{m_{k+1}})\oplus\cdots\oplus\Z/(q_{m_{r}}),T)
\]
consisting of those $(r-k)$-tuples $(t_{m_{k+1}},\dots t_{m_{r}})$ in
$T$ such that $t_{m_{i}}\ne 1$ for $k+1\le i\le r$.
Using the Cayley map as in \cite[Section 7]{ACG} we can
find a $W$-equivariant homeomorphism
\[
T\setminus \{1\}\cong \t.
\]
Here $\t$ denotes the Lie algebra of $T$. Using this identification and
the restriction of the map $\varphi_{\m}$, we obtain a surjective map
\[
\psi_{\m}:G/T\times_{W} (\t^{k}\times Q(P_{\m}(\pi),T))\to
\Hom(P_{\m}(\pi),G)\setminus S_{1}(P_{\m}(\pi),SU(2))
\]
Moreover, $\psi_{\m}$ is injective except where the action of $W$ on
$\t^{k}\times Q(P_{\m}(\pi),T)$ is not free. We need to consider two
cases.

\begin{list}{\labelitemi}{\leftmargin=2em}
\item Suppose first that
that $q_{m_{i}}$ is odd for some $k+1\le i\le r$.
In that case $W$ acts freely on $ \t^{k}\times Q(P_{\m}(\pi),T)$
and we have a $W$-equivariant homeomorphism
\[
\t^{k}\times Q(P_{\m}(\pi),T)\cong
\bigsqcup_{A(\m,\pi)}\t^{k}\times W
\]
Here $A(\m,\pi)$ is the number of $W$-orbits in $Q(P_{\m}(\pi),T)$.
This number is precisely
\[
A(\m,\pi)=\frac{1}{2}(q_{m_{k+1}}-1)\cdots(q_{m_{r}}-1).
\]
Therefore
\[
(G/T\times_{W} (\t^{k}\times Q(P_{\m}(\pi),T)))^{+}\cong
(\bigsqcup_{A(\m,\pi)}G/T\times_{W} (\t^{k}\times W))^{+}\cong
\bigvee_{A(\m,\pi)}(G/T\times\t^{k})^{+}.
\]
Note that $G/T=\S^{2}$ and thus
\[
(G/T\times \t^{k})^{+}\cong \Sigma^{k}(\S^{2}_{+})
\cong \S^{k+2}\vee \S^{k}.
\]
In this case the map $\psi_{\m}$ is a homeomorphism as the action of
$W$ on $ \t^{k}\times Q(P_{\m}(\pi),T)$ is free.
This shows that if $q_{m_{i}}$ is odd for some $k+1\le i\le r$
then
\[
\Hom(P_{\m}(\pi),SU(2))/S_{1}(P_{\m}(\pi),SU(2))\cong
(\bigvee_{A(\m,\pi)} \S^{k+2})\vee (\bigvee_{A(\m,\pi)} \S^{k}).
\]

\item Suppose now that
$P_{\m}(\pi)=\Z^{k}\oplus\Z/(q_{m_{k+1}})\cdots\oplus \Z/(q_{m_{r}})$
and that $q_{m_{i}}$ is even for every $k+1\le i\le r$. In this
case we have two kinds of elements in $Q(P_{\m}(\pi),T)$. On the one
hand we have the $(r-k)$-tuple $(-1,\dots,-1)$ on which $W$ acts
trivially. For all other elements in $Q(P_{\m}(\pi),T)$ the action of
$W$ is free. This shows that there is a $W$-equivariant homeomorphism
\[
\t^{k}\times Q(P_{\m}(\pi),T)\cong
\t^{k}\sqcup(\bigsqcup_{A(\m,\pi)}\t^{k}\times W).
\]
Here $A(\m,\pi)$ denotes the number of $W$-orbits in $Q(P_{\m}(\pi),T)$
different from the trivial orbit  $[(-1,\dots,-1)]$. This number is
precisely
\[
A(\m,\pi)=\frac{1}{2}\left((q_{m_{1}}-1)\cdots(q_{m_{r}}-1)-1\right).
\]
The map $\psi_{\m}$ is no longer injective. Note that $G/T\times_{W}\t^{k}$
is the Whitney sum of $k$-copies of the canonical vector bundle
over $\R P^{2}$ and $\psi$ maps the zero section onto the
$n$-tuple $(-1,\dots,-1)$. On the other hand, the restriction of
$\psi_{\m}$ onto the factor
\[
G/T\times_{W}\left(\bigsqcup_{A(\m,\pi)}\t^{k}\times W\right)
\cong \bigsqcup_{A(\m,\pi)}G/T\times\t^{k}
\]
is injective. This shows that if $q_{m_{i}}$ is even for every
 $k+1\le i\le r$ then
\[
\Hom(P_{\m}(\pi),SU(2))/S_{1}(P_{\m}(\pi),SU(2))\cong
(\R P^{2})^{k\lambda_{2}}/s_{k}(\R P^{2})\vee
(\bigvee_{A(\m,\pi)} \S^{k+2})\vee (\bigvee_{A(\m,\pi)} \S^{k}).
\]
\end{list}
\end{list}
\end{example}

\medskip

We can use the above for example to establish  the stable
homotopy type of the space of homomorphisms
$\Hom(\Z^{2}\oplus\Z/(2)\oplus \Z/(3),SU(2))$. In this
case we have that after one suspension
$\Hom(\Z^{2}\oplus\Z/(2)\oplus \Z/(3),SU(2))$
is homotopy equivalent to
\[
(\bigvee^{2}(\R P^{2})^{2\lambda_{2}}/s_{2}(\R P^{2}))
\vee(\bigvee^{2}\S^{4})\vee(\bigvee^{8}\S^{3})\vee
(\bigvee^{2}\S^{2})\vee(\bigvee^{2}\S_{+}^{2})\vee
(\bigvee^{4}\S^{1})\vee\S^{0}.
\]

\medskip

\begin{example}
Suppose that $\pi=\Z\oplus A$, where $A$ is a
finite abelian group.  Choose $G\in \P$ and assume that
$A$ is such that  the action of $W$ on
$\Hom(A,T)\setminus \{1\}$
is free. Since $W$ fixes the trivial homomorphism
$1\in \Hom(A,T)$,  then  the decomposition of
$\Hom(A,T)$ into $W$-orbits shows that in particular
$|W|$ divides $(|\Hom(A,T)|-1)$ under this assumption.

We will show that in this case
\[
\Sigma \Hom(\pi,G)\simeq \Sigma(\bigvee_{k} T)\vee
\Sigma(\bigvee_{k} G/T\wedge T)\vee \Sigma G
\vee (\bigsqcup_{k}G/T )_{+}
\]
Here $k:= (|\Hom(A,T)|-1)/|W|$ is the number of distinct $W$-orbits
on the set $\Hom(A,T)$ that are different from the
one corresponding to the trivial homomorphism.

Indeed,  using Theorem \ref{stable splitting} we obtain a
homotopy equivalence
\[
\Sigma \Hom(\pi,G)\simeq \Sigma \Hom(\pi,G)/S_{1}(\pi,G)\vee
\Sigma \Hom(\Z,G)/S_{1}(\Z,G)\vee
\Sigma \Hom(A,G)/S_{1}(A,G).
\]
Trivially $\Hom(\Z,G)/S_{1}(\Z,G)=G$. Also, since
$A$ is a finite abelian group then by Proposition
\ref{finite groups} we have
\[
\Hom(A,G)\cong \bigsqcup_{[f]\in \Hom(A,T)/W}G/G_{f}.
\]
Here $[f]$ runs through all the $W$-orbits in the
finite set $\Hom(A,T)$ and $G_{f}$ is a maximal rank
subgroup such that $W(G_{f})=W_{f}$. In $\Hom(A,T)$
we have two different kinds of orbits. On the one hand, we have
the orbit corresponding to the trivial homomorphism in
$\Hom(A,T)$. For this orbit we have $W_{f}=W$ and
$G_{f}=G$. The assumptions on $A$ imply that
for all other orbits in $\Hom(A,T)/W$
we have $W_{f}=1$ and thus $G_{f}=T$. This shows that
\[
\Hom(A,G)/S_{1}(A,G)=
(\bigsqcup_{k}G/T )_{+}.
\]
We now determine the stable factor $ \Hom(\pi,G)/S_{1}(\pi,G)$.
For this consider the map
\[
\varphi:G/T\times_{W} \Hom(\pi,T)\to \Hom(\pi,G).
\]
Since the action of $G$ on $\Hom(\pi,G)$ has maximal
rank isotropy subgroups $\varphi$ is surjective.  Moreover,
the restriction of $\varphi$
induces a surjective map
\[
\varphi_{|}: G/T\times_{W} ((T\setminus \{1\})\times (\Hom(A,T)\setminus \{1\}))
\to \Hom(\pi,G)\setminus S_{1}(\pi,G).
\]
Since the action of $W$ on $\Hom(A,T)\setminus \{1\}$ is
free we have that this restriction map is a homeomorphism.
Also
\[
G/T\times_{W} ((T\setminus \{1\})\times (\Hom(A,T)\setminus \{1\}))\cong
\bigsqcup_{k}G/T\times (T\setminus \{1\}).
\]
This shows that
\[
\Hom(\pi,G)/S_{1}(\pi,G)\cong \bigvee_{k} (G/T\times (T\setminus \{1\}))^{+}
\]
Note that
\[
(G/T\times (T\setminus \{1\}))^{+}\cong (G/T\times T)/(G/T\times \{1\}).
\]
and it is easy to see that there is a homotopy equivalence
\[
\Sigma ((G/T\times T)/(G/T\times \{1\}))\simeq \Sigma T\vee
\Sigma G/T\wedge T.
\]
This shows that
\[
\Sigma \Hom(\pi,G)/S_{1}(\pi,G)\cong \Sigma(\bigvee_{k} T)\vee
\Sigma(\bigvee_{k} G/T\wedge T)
\]
proving the claim.
\end{example}

\section{Fundamental group}\label{fundamental group}

In this section we study the fundamental group of the spaces of
homomorphisms $\Hom(\pi,G)$ under different choices of base point.

\medskip

Suppose first that $\pi=\Z^{n}$ and that $G$ is a compact Lie group.
Let $\BONE\in \Hom(\Z^{n},G)$ be the trivial representation. If we
give $\Hom(\Z^{n},G)$ the base point $\BONE$ then by \cite[Theorem
1.1]{GPS} there is a natural isomorphism
$\pi_{1}(\Hom(\Z^{n},G))\cong (\pi_{1}(G))^{n}$. The tools applied in 
\cite{GPS} can be used to generalize this result to the class of spaces
of homomorphisms $\Hom(\pi,G)$.  Here we  need to assume that $\pi$
is a finitely generated abelian group and that $G$ is a Lie group 
in the class $\P$. Under these assumptions \cite[Theorem 1.1]{GPS} can be
generalized for any choice of base point in $\Hom(\pi,G)$. Write
$\pi$ in the form $\pi=\Z^{n}\oplus A$, where $A$ is a finite
abelian group. Suppose first that $n=0$ and thus $\pi$ is a finite
group. In this case by Proposition  \ref{finite groups} there is a
homeomorphism
\[
\Phi:\Hom(\pi,G)\to \bigsqcup_{[f]\in \Hom(\pi,T)/W}G/G_{f},
\]
where each $G_{f}$ is a maximal rank isotropy subgroup
with $W(G_{f})=W_{f}$. For each maximal rank subgroup
$H\subset G$ we have $\pi_{1}(G/H)=1$. It follows that
$\pi_{1}(\Hom(\pi,G))=1$ for any choice of base point in this
case.  This handles the case of finite groups. Suppose then that
$n\ge 1$.   Let $T\subset G$ be a maximal torus.  Note that
$\Hom(\pi,G)^{T}=\Hom(\pi,T)$
and since $T$ is abelian we have
$\Hom(\pi,T)= \Hom(\Z^{n},T)\times \Hom(A,T)
=T^{n}\times \Hom(A,T).$
Choose $f\in \Hom(A,T)$ and let $\BONE\in \Hom(\Z^{n},T)$ denote
the trivial representation.  Let
\[
\BONE_{f}:=\BONE\times f\in \Hom(\Z^{n},T)\times \Hom(A,T)=
\Hom(\pi,T)\hookrightarrow \Hom(\Z^{n},G)
\]
and denote by $\Hom(\pi,G)_{\BONE_{f}}$ the path--connected component
of $\Hom(\pi,G)$ containing $\BONE_{f}$. It is easy to see that
\[
\Hom(\pi,G)=\bigsqcup_{[f]\in \Hom(A,T)/W}\Hom(\pi,G)_{\BONE_{f}}.
\]
Since the fundamental group does not depend, up to isomorphism,
on the base point chosen on a path--connected space, it
suffices to compute $\pi_{1}(\Hom(\pi,G)_{\BONE_{f}})$ for any
$f\in \Hom(A,T)$, where $\Hom(\pi,G)_{\BONE_{f}}$ is given the
base point $\BONE_{f}$.  Fix $f\in\Hom(A,T)$. Note that
$\Hom(\pi,G)_{\BONE_{f}}$ is invariant under the conjugation action of
$G$ and this action has connected maximal isotropy subgroups.
In this case
$\Hom(\pi,G)_{\BONE_{f}}^{T}\cong T^{n}\times W/W_{f}$.
As pointed out in Section \ref{max rank} we have a surjective
map
\[
\varphi:G/T\times_{W}\Hom(\pi,G)_{\BONE_{f}}^{T}\cong
G/T\times_{W_{f}}T^{n}\to \Hom(\pi,G)_{\BONE_{f}}
\]
that has connected fibers. As observed before this map
is not injective in general; however, there is a large set on
which this has this desirable property. Let
$\Fi(\pi,f)$  be the subset of $ \Hom(\pi,G)_{\BONE_{f}}^{T}$
on which $W$ acts freely.  Then the restriction of $\varphi$
\[
\varphi_{|}:G/T\times_{W}\Fi(\pi,f)\to \varphi(G/T\times_{W}\Fi(\pi,f))
\subset\Hom(\pi,G)_{\BONE_{f}}
\]
 is a homeomorphism onto its image.  Assume further
that $G$ is not a torus. Then under
this additional assumption the complement of
$G/T\times_{W}\Fi(\pi,f)$ is an analytic subspace
of $G/T\times_{W}(T^{n}\times W/W_{f})$ of
co-dimension at least $2$. Indeed, if $G$ is not
a torus then  $G/T$ is a smooth manifold of dimension
at least $2$ and $\Hom(\pi,G)_{\BONE_{f}}^{T}\setminus \Fi(\pi,f)$
submanifold of $\Hom(\pi,G)_{\BONE_{f}}^{T}$
of co-dimension at least $1$. This proves the claim.
Note that when $G$ is a torus then $W$ is trivial,
$\Fi(\pi,f)=\Hom(\pi,G)_{\BONE_{f}}^{T}$  and the map
$\varphi$ is a homeomorphism.

\medskip

Following \cite{GPS} we have the following definition.

\begin{definition}
Define $\H^{r}_{f}$ to be the image of $G/T\times_{W} \Fi(\pi,f)$ under
the map $\varphi$. We refer to $\H^{r}_{f}$ as the regular part
of $\Hom(\pi,G)_{\BONE_{f}}$. Also define
$\H^{s}_{f}:=\Hom(\pi,G)_{\BONE_{f}}\setminus \H^{r}_{f}$.
We refer to $\H^{s}_{f}$ as the singular part of
$\Hom(\pi,G)_{\BONE_{f}}$.
\end{definition}

Note that  $\Hom(\pi,G)_{\BONE_{f}}$ is a  compact real analytic space
and since $\H^{s}_{f}$ is the image of the compact analytic space
$(G/T\times_{W} T^{n}\times W/W_{f})\setminus (G/T\times_{W} \Fi(\pi,f))$
under the analytic map $\varphi$, it follows that $\H^{s}_{f}$ is a
compact analytic subspace of $\Hom(\pi,G)_{\BONE_{f}}$.
As a consequence of the Whitney stratification
theorem it follows that $\Hom(\pi,G)_{\BONE_{f}}$ can be given the
structure of a simplicial complex in such a way that $\H^{s}_{f}$ is a
subcomplex.  Note that when $G$ is a torus $\H_{f}^{s}$ is empty.
On the other hand, if $G$ is not a torus then using the fact that
the complement of  $G/T\times_{W}\Fi(\pi,f)$ is an analytic subspace of
$G/T\times_{W}(T^{n}\times W/W_{f})$ of  co-dimension at least
$2$ and an argument similar to the one provided in
\cite[Lemma 2.4]{GPS} the following lemma is obtained for
any $G\in \P$.

\begin{lemma}\label{nowhere dense}
The space $\Hom(\pi,G)_{\BONE_{f}}$ is a compact simplicial
complex and the singular part $\H^{s}_{f}$ is a subcomplex.
Also, $\H^{s}_{f}$ is nowhere dense and does not
disconnect connected open subsets of $\Hom(\pi,G)_{\BONE_{f}}$.
\end{lemma}

The previous lemma can be used as a first step for the
computation of $\pi_{1}(\Hom(\pi,G)_{\BONE_{f}})$. Indeed,
suppose that $X$ is a compact simplicial complex and that
$Y\subset X$ is a subcomplex. Assume that $X\setminus Y$ is
dense and that $Y$ does not separate any connected open set in $X$.
If $x_{0}\in X\setminus Y$ is the base point  then by
\cite[Lemma 2.5]{GPS} the inclusion map $i:X\setminus Y\to X$
induces a surjective homomorphism
$i_{*}:\pi_{1}(X\setminus Y,x_{0})\to \pi_{1}(X,x_{0})$.
This can be applied in our situation.  Choose
$x_{0}\in \H_{f}^{r}$ as the base point.  Using
Lemma \ref{nowhere dense} we obtain that
the inclusion
$i:\H^{r}_{f}\hookrightarrow \ \Hom(\pi,G)_{\BONE_{f}}$
induces a surjective homomorphism
$i_{*}:\pi_{1}(\H^{r}_{f})\to \pi_{1}(\Hom(\pi,G)_{\BONE_{f}})$.
The same argument shows that the inclusion map
\[
i:G/T\times_{W} \Fi(\pi,f)\hookrightarrow
G/T\times_{W} \Hom(\pi,G)_{\BONE_{f}}^{T}
\]
induces a surjective homomorphism
\[
i_{*}:\pi_{1}(G/T\times_{W} \Fi(\pi,f))\to
\pi_{1}(G/T\times_{W} \Hom(\pi,G)_{\BONE_{f}}^{T}).
\]
Since $\varphi_{|}:G/T\times_{W} \Fi(\pi,f)\to \H^{r}$ is a
homeomorphism and the fundamental group of a connected space does
not depend on the choice of base point, up homeomorphism,  we obtain
the following proposition (compare \cite[Corollary 2.6]{GPS}).

\begin{proposition}\label{surjective}
Suppose that $G\in \P$ and that $\pi$ is a finitely generated
abelian group. Then the map
\[
\varphi: G/T\times_{W} \Hom(\pi,G)_{\BONE_{f}}^{T}
\to \Hom(\pi,G)_{\BONE_{f}}
\]
is $\pi_{1}$-surjective.
\end{proposition}

Note that
\[
G/T\times_{W} \Hom(\pi,G)_{\BONE_{f}}^{T}\cong
G/T\times_{W_{f}}T^{n}.
\]
Since $W_{f}$ acts freely on $G/T$ the projection map
$p$ induces a fibration sequence
\[
T^{n}\to G/T\times_{W_{f}} T^{n}\stackrel{p}{\rightarrow}
(G/T)/W_{f}\cong G/N_{G_{f}}(T).
\]
The tail of the homotopy long exact sequence
associated  to this fibration is the exact sequence
\begin{equation}\label{exact sequence}
\pi_{1}(T^{n})\to \pi_{1}( G/T\times_{W_{f}} T^{n})
\to \pi_{1}(G/N_{G_{f}}(T))\to 1.
\end{equation}
Note that $\BONE\in T^{n}$ is a fixed point of $W_{f}$. Therefore
the map
\begin{align*}
s:G/N_{G_{f}}(T)&\to G/T\times_{W_{f}} T^{n}\\
[g]&\mapsto [g\times \BONE]
\end{align*}
is a section of $p$ and in particular the sequence (\ref{exact sequence})
splits. This proves that $\pi_{1}(G/T\times_{W_{f}} T^{n})$ is
generated by $\pi_{1}(T^{n})$ and $s_{*}(\pi_{1}(G/N_{G_{f}}(T)))$.
Next we prove the following lemma.

\begin{lemma}
If $\alpha:[0,1]\to G/N_{G_{f}}(T)$ is a loop then
$\varphi\circ s\circ \alpha$ is homotopic to the the trivial
loop in $\Hom(\pi,G)_{\BONE_{f}}$. Therefore
$s_{*}(\pi_{1}(G/N_{G_{f}}(T)))\subset \Ker(\varphi_{*})$.
\end{lemma}
\Proof
Let $\alpha:[0,1]\to G/N_{G_{f}}(T)$ be a loop. Note that
$\Hom(\pi,G)_{\BONE_{f}}$ can be seen
as a subspace of $\Hom(\Z^{n},G)\times \Hom(A,G)$.
Under this identification  $\beta:=\varphi\circ s\circ \alpha$
is the loop in $ \Hom(\pi,G)_{\BONE_{f}}$ given by
\begin{align*}
\beta:=\varphi\circ s\circ \alpha:[0,1]&\to \Hom(\pi,G)_{\BONE_{f}}
\subset \Hom(\Z^{n},G)\times \Hom(A,G)\\
t&\mapsto (\BONE, \alpha(t)f\alpha(t)^{-1}).
\end{align*}
Let $G_{f}=Z_{G}(f)$ be the  subspace of elements in $G$
commuting with $f(x)$ for all $x\in A$ and $G\cdot f$
the space of elements in $\Hom(A,G)$ conjugated to $f$.
Then
\[
\beta:[0,1]\to \{\BONE\}\times G\cdot f\subset
\Hom(\pi,G)_{\BONE_{f}}.
\]
There is a homeomorphism $G\cdot f\cong G/G_{f}$ and
$G_{f}$ is a maximal rank subgroup in $G$ as $T\subset G_{f}$.
In particular the homogeneous space $G/G_{f}$ is simply connected.
The simply connectedness of  $G\cdot f $ shows that
up to homotopy $\beta$ is the  trivial loop
in $\Hom(\pi,G)_{\BONE_{f}}$ proving the lemma.
\qed

\medskip

The previous lemma together with Proposition \ref{surjective}
and the fact that  $\pi_{1}(G/T\times_{W_{f}} T^{n})$ is
generated by $\pi_{1}(T^{n})$ and $s_{*}(\pi_{1}(G/N_{G_{f}}(T)))$
show that the map
\begin{align*}
\sigma_{f}:T^{n}&\to \Hom(\pi,G)_{f}\subset
\Hom(\Z^{n},G)\times \Hom(A,G)\\
(t_{1},\dots,t_{n})&\mapsto(t_{1},\dots,t_{n})\times \{f\}
\end{align*}
is $\pi_{1}$-surjective. On the other hand,
the inclusion $T\subset G_{f}$  shows that
$T^{n}\subset \Hom(\Z^{n},G_{f})$ and there is a commutative
diagram
\begin{equation}\label{diagram}
\xymatrix{
T^{n}\ar[r]^{\sigma_{f}}\ar[d]   &\Hom(\pi,G)_{\BONE_{f}}.\\
\Hom(\Z^{n},G_{f})\ar[ru]_{i_{f}}        &
}
\end{equation}
In the previous commutative diagram $i_{f}$ denotes the
map
\begin{align*}
i_{f}:\Hom(\Z^{n},G_{f})&\to \Hom(\pi,G)_{f}\subset
\Hom(\Z^{n},G)\times \Hom(A,G)\\
(x_{1},\dots x_{n})&\mapsto (x_{1},\dots x_{n})\times\{f\}.
\end{align*}

The inclusion map $T\subset G_{f}$ is $\pi_{1}$-surjective
and by \cite[Theorem 1.1]{GPS} the map induced by the inclusion
$\pi_{1}(\Hom(\Z^{n},G_{f}))\to \pi_{1}(G_{f}^{n})\cong (\pi_{1}(G_{f}))^{n}$
is an isomorphism. This proves that the map
$\pi_{1}(T^{n})\to \pi_{1}(\Hom(\Z^{n},G_{f}))$ is surjective.
Using the commutativity of diagram (\ref{diagram}) and the
fact that $\sigma_{f}$ is $\pi_{1}$-surjective we obtain the
following corollary.

\begin{corollary}\label{inclusion surjective}
Suppose that $\pi$ is a finitely generated abelian . Then the map
\[
i_{f}:\Hom(\Z^{n},G_{f})\to \Hom(\pi,G)_{\BONE_{f}}
\]
described above is $\pi_{1}$-surjective.
\end{corollary}

We are now ready to prove the following theorem which is the
main theorem of this section.

\begin{theorem}
Let $\pi=\Z^{n}\oplus A$, with $A$ a finite abelian group
and $G\in\P$. Let $f\in \Hom(A,T)$ and let
$\BONE:=\BONE\times f\in \Hom(\pi,G)$ be the base point of
$ \Hom(\pi,G)$. Then there is a natural
isomorphism
$\pi_{1}(\Hom(\pi,G))\cong (\pi_{1}(G_{f}))^{n}$,
where $G_{f}=Z_{G}(f)$ is the subgroup of elements in $G$ commuting
with $f(x)$ for all $x\in A$.
\end{theorem}
\Proof
Suppose first that $\pi$ is a finite group and thus $n=0$.
Then as proved above $\pi_{1}(\Hom(\pi,G))=1$ for any choice of
base point and the theorem is true in this case. Suppose then
that $n\ge 1$. By Corollary \ref{inclusion surjective} the map
$i_{f}$ is $\pi_{1}$-surjective. We now show that in fact
$(i_{f})_{*}:\pi_{1}(\Hom(\Z^{n},G_{f}))\to
\pi_{1}( \Hom(\pi,G)_{\BONE_{f}})$
is an isomorphism. This together with
the isomorphism
$\pi_{1}(\Hom(\Z^{n},G_{f}))\cong (\pi_{1}(G_{f}))^{n}$
provided by \cite[Theorem 1.1]{GPS} proves the theorem.

To start note that $G_{f}$ is such that $\pi_{1}(G_{f})$ is torsion
free. Therefore we can write $\pi_{1}(G_{f})=\Z^{a}$ for
some integer $a$ and the map
\[
(i_{f})_{*}:\pi_{1}(\Hom(\Z^{n},G_{f}))\cong \Z^{na}\to
\pi_{1}( \Hom(\pi,G)_{\BONE_{f}})
\]
is a surjection. This shows in particular that
$\pi_{1}( \Hom(\pi,G)_{\BONE_{f}})$ is abelian
and of rank at most $na$. We are going to show that
in fact
\[
r:={\rm{rank}}_{\Z}(\pi_{1}( \Hom(\pi,G)_{\BONE_{f}}))=na.
\]
The only way this is possible is that $(i_{f})_{*}$ is an
isomorphism, proving the theorem. We now verify this.
The universal coefficient theorem together
with the Hurewicz theorem provide an isomorphism
of $\Q$-vector spaces
\begin{align}\label{ranks}
H^{1}(\Hom(\pi,G)_{\BONE_{f}};\Q)&\cong \Q^{r},\\
H^{1}(\Hom(\Z^{n},G_{f});\Q)&\cong \Q^{na}.
\end{align}
On the other hand, since the conjugation action of
$G$ on $\Hom(\pi,G)_{\BONE_{f}}$ has connected
maximal rank isotropy subgroups, then by Theorem
\ref{rational cohomology} there is an isomorphism
\[
H^{*}(\Hom(\pi,G)_{\BONE_{f}};\Q)\cong
H^{*}(G/T\times_{W}(\Hom(\pi,G)_{\BONE_{f}})^{T};\Q).
\]
In this case
\[
G/T\times_{W}(\Hom(\pi,G)_{\BONE_{f}})^{T}\cong
G/T\times_{W_{f}}T^{n}.
\]
In particular
\begin{equation}\label{iso1}
H^{1}(\Hom(\pi,G)_{\BONE_{f}};\Q)\cong
H^{1}(G/T\times T^{n};\Q)^{W_{f}}\cong
H^{1}(T^{n};\Q)^{W_{f}}.
\end{equation}
The second isomorphism follows from the
fact that $H^{1}(G/T;\Q)=0$ as $G/T$ is simply
connected. On the other hand,  by \cite[Theorem 1.1]{Hauschild}
it follows $W(G_{f})=W_{f}$. The conjugation action of
$G_{f}$ on $\Hom(\Z^{n},G_{f})$ also has maximal rank
isotropy subgroups. The same argument as above yields the
following isomorphisms
\begin{equation}\label{iso2}
H^{1}(\Hom(\Z^{n},G_{f});\Q) \cong
H^{1}(G_{f}/T\times T^{n};\Q)^{W_{f}}\cong
H^{1}(T^{n};\Q)^{W_{f}}.
\end{equation}
Equations (\ref{iso1}) and (\ref{iso2}) show that there is an
isomorphism of $\Q$-vector spaces
\[
H^{1}(\Hom(\Z^{n},G_{f});\Q)\cong
 H^{1}(\Hom(\pi,G)_{\BONE_{f}};\Q)
\]
and thus $n=ra$ by (\ref{ranks}).
\qed

\section{Equivariant $K$-theory}\label{equivariant K-theory}

In this section we study the $G$-equivariant $K$-theory of the spaces
of homomorphisms $\Hom(\pi,G)$. We divide our study according to
the nature of the group $\pi$.  From now on we fix $T$ a maximal
torus in $G$ and let $W$ be the associated Weyl group.

\subsection{Finite abelian groups}

We first consider the case where $\pi$ is a finite abelian group.

\medskip

Fix a finite abelian group $\pi$ and $G\in \P$. By Proposition 
\ref{finite groups} there is a $G$-equivariant homeomorphism
\[
\Phi:\Hom(\pi,G)\to \bigsqcup_{[f]\in \Hom(\pi,T)/W}G/G_{f}.
\]
Given a subgroup $H\subset G$ we have
\begin{equation*}
K^{q}_{G}(G/H)\cong \left\{
\begin{array}{ccc}
R(H)&\text{ if } &q \text{ is even},\\
0& \text{ if } &q \text{ is odd}.
\end{array}
\right.
\end{equation*}
By \cite[Theorem 1]{Pittie} it follows that if $H\subset G$ is a 
subgroup of maximal rank then $R(H)$ is a free module over $R(G)$ 
of rank $|W|/|WH|$. As a corollary of this the following is obtained.

\begin{corollary}\label{rank 0}
Let $G\in \P$ and $\pi$ be a finite abelian group.
Then $K_{G}^{0}(\Hom(\pi,G))$ is a free $R(G)$-module
of rank $|\Hom(\pi,T)|$ and $K^{1}_{G}(\Hom(\pi,G))=0$.
\end{corollary}
\Proof
Using Proposition \ref{finite groups} and the above we have
$K_{G}^{1}(\Hom(\pi,G))=0$ and
\[
K^{0}_{G}(\Hom(\pi,G))\cong \bigoplus_{[f]\in \Hom(\pi,T)/W}R(G_{f}).
\]
Each $R(G_{f})$ is a free module over $R(G)$ of rank
$|W|/|W_{f}|$. Note that  $W(G_{f})=W_{f}$, where $W_{f}$
denotes the isotropy subgroup at $f$, under the action of $W$ on the 
finite set $\Hom(\pi,T)$. The partition of $\Hom(\pi,T)$ into 
the different $W$-orbits provides the identity
\[
|\Hom(\pi,T)|=\sum_{[f]\in \Hom(\pi,T)/W} |W|/|W_{f}|.
\]
This proves that $K^{0}_{G}(\Hom(\pi,G))$ is free as a module
over $R(G)$ of rank $|\Hom(\pi,T)|$.
\qed

\medskip

\noindent{\bf{Remark:}} The previous corollary is not true in general 
if $G\notin \P$. For example, it can be seen that
$K^{*}_{PU(3)}(\Hom((\Z/(3))^{2},PU(3)))$ is not free as a module over
$R(PU(3))$.

\subsection{Abelian groups of rank one}

We now consider the case where $\pi$ is a finitely
generated abelian group of rank one. Thus we can write $\pi$
in the form $\pi=\Z\oplus A$ where $A$ is a
finite abelian group.

\medskip

Suppose that  $X$ is a $G$-CW complex. The skeleton filtration of $X$
induces a multiplicative spectral sequence  (see \cite{Segal}) with
\begin{equation}\label{spectral sequence}
E_{2}^{p,q}=H_{G}^{p}(X;\K^{q}_{G})\Longrightarrow K^{p+q}_{G}(X).
\end{equation}
The $E_{2}$-term of this spectral sequence is the Bredon cohomology
of $X$ with respect to the coefficient system $\K_{G}^{q}$
defined by $G/H\mapsto K^{q}_{G}(G/H)$.

Suppose that $G\in \P$ and that $\pi=\Z\oplus A$, where $A$
is a finite abelian group. Then Corollary \ref{CW complex} gives
$X:=\Hom(\pi,G)$ the structure of a $G$-CW complex and we can use
the previous spectral sequence to compute $K_{G}^{*}(\Hom(\pi,G))$.
In \cite[Theorem 1.6]{AG} a criterion for the collapse of the spectral 
sequence (\ref{spectral sequence}) without extension problems
was provided. This criterion can be used in this case to compute
the structure of $K_{G}^{*}(\Hom(\pi,G))$ as a module over $R(G)$.
Let $\Phi$ be the root system associated to
$(G,T)$.  Fix a subset $\Phi^{+}$ of positive roots
of $\Phi$ and let $\Delta=\{\alpha_{1},\dots,\alpha_{r}\}$ be an
ordering of the corresponding set of simple roots.
Suppose that $W_{i}\subset W$ is a reflection subgroup.
Let $\Phi_{i}$ be the corresponding root system and $\Phi_{i}^{+}$ the
corresponding positive roots. Define
\[
W_{i}^{\ell}:=\{w\in W ~|~  w(\Phi^{+}_{i})\subset \Phi^{+}\}.
\]
The set $W_{i}^{\ell}$ forms a system of representatives of the left
cosets in $W/W_{i}$ by \cite[Lemma 2.5]{Steinberg}. In a precise
way, this means that any element $w\in W$ can be factored in a
unique way in the form $w= ux \text{ with } u\in W_{i}^{\ell} \text{
and } x\in W_{i}$. In order to apply  the criterion provided in
\cite[Theorem 1.6]{AG} we must verify the hypothesis required there.
In particular we need to show that $X^{T}$ has the structure of a
$W$-CW complex in such a way that there is a CW-subcomplex $K$ of
$X^{T}$ such that for every element $x\in X^{T}$ there is a unique 
$w\in W$ such that $wx\in K$. We construct such a subcomplex next.
To start note that $X^{T}=\Hom(\pi,G)^{T}=T\times \Hom(A,T)$ and
$\Hom(A,T)$ is a discrete set endowed with an action of $W$.  If we
assume  that $G$ is simply connected then as pointed out in
\cite[Section 6.1]{AG} the (closed) alcoves in $T$ provide a
structure of a $W$-CW complex in $T$ in such a way that $K(\Delta)$,
the alcove determined by $\Delta$, is a sub CW-complex of $T$
such that any element in $T$ has a unique representative in 
$K(\Delta)$ under the $W$-action. Moreover, for each cell 
$\sigma$ in $K(\Delta)$, the isotropy subgroup $W_{\sigma}$ 
is a reflection subgroup of the form $W_{I}$ for some 
$I\subset \Delta$. Here $W_{I}$ denotes the
reflection subgroup generated by the reflections of the form
$s_{\alpha}$  for $\alpha \in I$. This can be used  to produce a sub
CW-complex in $\Hom(\pi,G)^{T}$ satisfying similar properties in the
following way. Let $f_{1},\dots, f_{m}$ be a set of representatives
for the $W$-orbits in the  discrete set $\Hom(A,T)$. We can choose
each $f_{i}$ in such a way that the isotropy subgroup $W_{f_{i}}$ is
a reflection subgroup of $W$ of the form $W_{I_{i}}$ for some
$I_{i}\subset \Delta$. For each $1\le i\le m$ let $W_{f_{i}}^{\ell}$
be a system of minimal length representatives of $W/W_{f_{i}}$ as
defined above. Define
\[
L(\Delta):=\bigcup_{i=1}^{m}\bigcup_{u\in W_{f_{i}}^{\ell}}
\left(u^{-1}K(\Delta)\times \{f_{i}\}\right)\subset T\times \Hom(A,T)=X^{T}.
\]
Defined in this way $L(\Delta)\subset X^{T}$ is a sub CW-complex. We
now show that $L(\Delta)$ is such that for every element $x\in X^{T}$ there 
is a unique $w\in W$ such that $wx\in L(\Delta)$. To see this, since
$K(\Delta)\subset T$ satisfies this property, it suffices to see that 
for any $i$ and any $v_{1},v_{2}\in W$ there are unique $v\in W$ and $u\in
W_{f_{i}}^{\ell}$ such that $v_{1}K(\Delta)\times
v_{2}f_{i}=v(u^{-1}K(\Delta)\times f_{i})$. Indeed, suppose that
$v_{1},v_{2}\in W$. Using the defining property of
$W_{f_{i}}^{\ell}$  we can find unique $u\in W_{f_{i}}^{\ell} $ and
$x\in W_{f_{i}}$ such that $v^{-1}_{1}v_{2}=ux$. Let $v=v_{1}u$.
Then $v_{1}=vu^{-1}$ and in particular
$v_{1}K(\Delta)=vu^{-1}K(\Delta)$. Also, $x=v^{-1}v_{2}\in
W_{f_{i}}$ and thus $v_{2}f_{i}=vf_{i}$. This shows that $v\in W$
and $u\in W_{f_{i}}^{\ell}$  are the unique elements such that
$v_{1}K(\Delta)\times v_{2}f_{i}=v(u^{-1}K(\Delta)\times f_{i})$. On
the other hand, note that $H^{*}(X^{T};\Z)$ is torsion--free and of
rank $2^{r}\cdot |\Hom(A,T)|$, where $r$ denotes the rank of the Lie
group $G$.  Also note that the isotropy subgroups of the action of
$W$ on $\Hom(\pi,G)^{T}=\Hom(\pi,T)$ are of the form $W_{I}$, with
$I\subset \Delta$.

The above work shows that the conditions of \cite[Theorem 1.6]{AG}
are satisfied yielding the next theorem.

\begin{theorem}\label{rank 1}
Suppose that $G\in \P$ is simply connected and of rank $r$.
Let $\pi=\Z\oplus A$ where $A$ is a finite abelian group.
Then $K_{G}^{*}(\Hom(\pi,G))$ is a free $R(G)$-module of rank
$2^{r}\cdot |\Hom(A,T)|$.
\end{theorem}

\noindent{\bf{Remark:}} Combining Corollary \ref{rank 0} and
Theorem \ref{rank 1} it follows that $K_{G}^{*}(\Hom(\pi,G))$ is free as a
module over $R(G)$ whenever $\pi$ is a finitely generated abelian
group of rank less or equal to $1$ and $G\in \P$ is simply connected.
As already pointed out in \cite{AG} this result does not extend to
all finitely generated abelian groups $\pi$ as
$K_{SU(2)}^{*}(\Hom(\Z^{2},SU(2)))$ contains torsion as a
$R(SU(2))$-module.

However, if we tensor with the rational numbers the previous
result does extend to the family of finitely generated abelian
groups  and all Lie groups $G\in \P$. This is done next.

\subsection{Finitely generated abelian groups}

We show that $K_{G}^{*}(\Hom(\pi,G))\otimes \Q$ is a free
$R(G)\otimes \Q$-module for all finitely generated abelian groups
$\pi$ and all Lie groups $G\in \P$.

\medskip

Let $G$ be a compact Lie group with $\pi_{1}(G)$ torsion--free
act on a compact space  $X$ with connected maximal rank
isotropy.  If we further assume that $X^{T}$ has the homotopy type of 
a $W$-CW complex then by  \cite[Theorem 1.1]{AG} 
$K^{*}_{G}(X)\otimes \Q$ is a free module 
over $R(G)\otimes \Q$ of rank equal to
$\sum_{i\ge 0}{\rm{rank}}_{\Q}H^{i}(X^{T};\Q)$. This theorem 
can be applied in our situation.  Let $\pi$ be a finitely generated
abelian group and $G\in \P$. Let $X:=\Hom(\pi,G)$. Then
the conjugation action of $G$ on $X$ has connected maximal
rank isotropy subgroups and $X$ has the homotopy type of a
$G$-CW complex by Proposition \ref{max rank for spaces of hom}
and Corollary \ref{CW complex}.  Note that 
$H^{*}(\Hom(\pi,G)^{T};\Q)$ is a $\Q$-vector space of
rank $2^{nr}\cdot |\Hom(A,T)|$, where $r$ is the rank of $G$. This
proves that the hypotheses of \cite[Theorem 1.1]{AG} are satisfied
in this case yielding the following.

\begin{theorem}
Suppose that $G\in \P$ is of rank $r$ and that $\pi$ is a finitely
generated abelian group written in the form $\pi=\Z^{n}\oplus A$,
where $A$ is a finite abelian group. Then
$K_{G}^{*}(\Hom(\pi,G))\otimes \Q$ is a free module over
$R(G)\otimes \Q$ of rank $2^{nr}\cdot |\Hom(A,T)|$.
\end{theorem}

\end{document}